\documentclass[12pt]{amsart}
\usepackage{amssymb,stmaryrd,color,xypic,mathdots,graphicx}

\theoremstyle{theorem}
\newtheorem{theorem}{Theorem}[section]
\newtheorem{lemma}[theorem]{Lemma}
\newtheorem{proposition}[theorem]{Proposition}
\newtheorem{corollary}[theorem]{Corollary}

\theoremstyle{definition}

\newtheorem{remark}[theorem]{Remark}

\newtheorem{question}[theorem]{Question}

\title[A simply connected universal fibration]{\large A \lowercase{simply connected universal fibration with unique path lifting over a} P\lowercase{eano continuum with non-simply connected universal covering space}}

\author{Jeremy Brazas}

\address{Department of Mathematics, West Chester University of Pennsylvania, West Chester, PA 19383, USA}

\email{jbrazas@wcupa.edu}

\author{Hanspeter Fischer}

\address{Department of Mathematical Sciences, Ball State University, Muncie, IN 47306, USA}

\email{hfischer@bsu.edu}

\date{May 10, 2024. \\ \mbox{\hspace{5pt} } 2020 {\em Mathematics Subject Classification.} Primary 55R05; Secondary 57M10}

\begin{document}

\begin{abstract}
We present a 2-dimensional Peano continuum $\mathbb{T}\subseteq \mathbb{R}^3$ with the following properties:
(1) There is a universal covering projection $q:\overline{\mathbb{T}}\rightarrow \mathbb{T}$ with uncountable fundamental group $\pi_1(\overline{\mathbb{T}})$;
(2) For every $1\not=[\overline{\alpha}]\in \pi_1(\overline{\mathbb{T}},\ast)$,
there is a covering projection $r:(E,e)\rightarrow (\overline{\mathbb{T}},\ast)$ such that $[\overline{\alpha}]\not\in r_\#\pi_1(E,e)$;
(3) There is no universal covering projection $r:E\rightarrow \overline{\mathbb{T}}$;
(4) The universal object $p:\widetilde{\mathbb{T}}\rightarrow \mathbb{T}$ in the category of fibrations with unique path lifting (and path-connected total space) over $\mathbb{T}$ has trivial fundamental group $\pi_1(\widetilde{\mathbb{T}})=1$;
(5) $p:\widetilde{\mathbb{T}}\rightarrow \mathbb{T}$ is not a path component of an inverse limit of covering projections over $\mathbb{T}$.
\end{abstract}

\maketitle

\section{introduction}
\renewcommand{\theequation}{$\mathcal S$}

\noindent Many foundational results of covering space theory can be developed in the more general context of (Hurewicz) fibrations with unique path lifting.
This is done, for example, in Spanier's classical textbook~\cite{S}, where in \S2.5 one also finds a sketch of the proof of the existence of a universal object $p:\widetilde{X}\rightarrow X$ in the category of fibrations with unique path lifting (and path-connected total space) over a fixed  path-connected base space $X$ and a proof of the following inclusion, assuming $X$ is also locally path connected: \begin{equation} \label{contained} p_\#\pi_1(\widetilde{X},\widetilde{x})\leqslant \bigcap_{{\mathcal U}\in Cov(X)} \pi({\mathcal U},x)\end{equation}
Here, $Cov(X)$ denotes the collection of all open covers of $X$ and $\pi({\mathcal U},x)$ denotes the subgroup of the fundamental group $\pi_1(X,x)$ generated by all elements of the form $[\gamma\cdot \delta\cdot \gamma^-]$ for some path $\gamma:([0,1],0)\rightarrow (X,x)$, its reverse $\gamma^-(t)=\gamma(1-t)$, and some loop $\delta:([0,1],\{0,1\})\rightarrow (U,\gamma(1))$ in some $U\in {\mathcal U}$.
We will denote this intersection by \[\pi^s(X,x):=\bigcap_{{\mathcal U}\in Cov(X)} \pi({\mathcal U},x)\] and refer to it as the {\em Spanier group} of $X$ at $x$. Since $p_\#$ is injective,
$\widetilde{X}$ is simply connected if $\pi^s(X,x)$ is trivial.

Given a connected and locally path-connected space $X$, we recall from \cite[2.5.11 \& 2.5.13]{S} that for a loop $\alpha:([0,1],\{0,1\})\rightarrow (X,x)$, we have that $[\alpha]\in \pi^s(X,x)$ if and only if for every covering projection $q:(E,e)\rightarrow (X,x)$, the unique lift $\alpha':([0,1],0)\rightarrow (E,e)$ of $\alpha$ with $q\circ \alpha'=\alpha$ is a loop, i.e. $\alpha'(1)=e$.

This raises the following natural question.

\begin{question}\label{question} Does equality hold in (\ref{contained})?
\end{question}

Since fibrations have continuous lifting of paths, it seems reasonable to speculate that the answer to Question~\ref{question} is ``yes'', as was recently claimed in \cite{CP}. However, we will show that this is not the case, in general, by presenting a Peano continuum $\mathbb{T}$ for which $\pi^s(\mathbb{T},\ast)$ is uncountable, although
the universal object $p:\widetilde{\mathbb{T}}\rightarrow \mathbb{T}$ in the category of fibrations with unique path lifting (and path-connected total space) over $\mathbb{T}$  satisfies
 $p_\#\pi_1(\widetilde{\mathbb{T}},\ast)=1$; in fact, $\mathbb{T}$ admits a (categorical) universal covering projection $q:(\overline{\mathbb{T}},\ast)\rightarrow (\mathbb{T},\ast)$ with $q_\#\pi_1(\overline{\mathbb{T}},\ast)=\pi^s(\mathbb{T},\ast)$ and $\pi^s(\overline{\mathbb{T}},\ast)=1$, while no universal covering projection over $\overline{\mathbb{T}}$ exists.

To put this result into perspective, we mention that
the Spanier group $\pi^s(X,x)$ is always contained in the kernel of the canonical homomorphism $\pi_1(X,x)\rightarrow \check{\pi}_1(X,x)$ to the first \v{C}ech homotopy group~\cite{FZ2007}, which is injective for spaces that are one-dimensional~\cite{EK} or planar~\cite{FZ2005},
and it equals this kernel when $X$ is locally path connected and metrizable \cite{BFa}. Therefore, $\pi^s(X,x)=1$ if $X$ is either one-dimensional or planar. Our Peano continuum  is a 2-dimensional subset of~$\mathbb{R}^3$. Also, if $X$ is
 locally path connected, semilocally simply connected, and metrizable, then $\pi_1(X,x)\rightarrow \check{\pi}_1(X,x)$ is an isomorphism \cite{FZ2007}. Our Peano continuum contains a $\pi_1$-injectively embedded shrinking wedge of countably many circles.

In contrast, the universal fibration with unique path lifting $p:\widetilde{X}\rightarrow X$ over the Griffiths space $X$ satisfies $p_\#\pi_1(\widetilde{X},\widetilde{x})=\pi^s(X,x)=\pi_1(X,x)$ with $\widetilde{X}=X$ and $p=id_X$ \cite[2.5.18]{S}, and $H_1(X)= \mathbb{Z}^\mathbb{N}/\bigoplus_\mathbb{N} \mathbb{Z}$~\cite{EF}.

For a comparison with other universal simply connected lifting objects, see Remark~\ref{GenUnivCov}.

\section{The space $\mathbb{T}$ and its universal covering $q:\overline{\mathbb{T}}\rightarrow \mathbb{T}$}

\noindent For each $m\in \mathbb{N}$, let $C_m\subseteq \mathbb{R}^2$ be the circle of radius $\frac{1}{m}$ centered at $(\frac{1}{m},0)$
and consider the ``infinite earring" $\mathbb{E}=\bigcup_{m\in \mathbb{N}} C_m\subseteq \mathbb{R}^2$ with basepoint $\ast=(0,0)\in \mathbb{E}$.
Let $\ell_m:[0,1]\rightarrow C_m$ be the loop at $\ast$ defined by
\[\ell_m(t)=(\frac{1}{m}(1-\cos 2\pi t), \frac{1}{m}\sin 2\pi t))\]
which traverses $C_m$ once clockwise.
Let $f:\mathbb{E}\rightarrow \mathbb{E}$ be the unique ``circle shifting'' map with \[f\circ \ell_m=\ell_{m+1}\] for all $m\in \mathbb{N}$.
 Let $\mathbb{T}$ be the mapping torus of $f:\mathbb{E}\rightarrow \mathbb{E}$. That is, let $\mathbb{T}$ be the quotient space $\mathbb{E}\times [0,1]/\sim$ with $(x,1)\sim (f(x),0)$ for all $x \in \mathbb{E}$.
 Let $h: \mathbb{E}\times [0,1]\rightarrow \mathbb{T}$ denote the quotient map. We will identify $\mathbb{E}$ with $h(\mathbb{E}\times \{0\}) \subseteq \mathbb{T}$ and let $\iota:(\mathbb{E},\ast)\rightarrow (\mathbb{T},\ast)$ denote inclusion.  Let $\eta:[0,1]\rightarrow \mathbb{E}\times [0,1]$ be given by $\eta(t)=(\ast,t)$ and put $\lambda=h\circ \eta$. Then $[\lambda]\in \pi_1(\mathbb{T},\ast)$.

\begin{lemma} The space $\mathbb{T}$ is a Peano continuum.
\end{lemma}

\begin{proof}
Note that $\mathbb{T}$ can be realized in $\mathbb{R}^3$ as the quotient of the Peano continuum $\mathbb{E}\times [0,1]$.
\end{proof}

For each $i\in \mathbb{Z}$, let $X_i$ be a copy of $\mathbb{E}$ and let $f_i=f:X_i\rightarrow X_{i+1}$.   Let $\overline{\mathbb{T}}$ be the mapping telescope of the bi-infinite sequence
 \[\cdots \stackrel{f_{-2}}{\longrightarrow} X_{-1} \stackrel{f_{-1}}{\longrightarrow} X_{0} \stackrel{f_{0}}{\longrightarrow} X_{1} \stackrel{f_{1}}{\longrightarrow}  X_{2} \stackrel{f_{2}}{\longrightarrow} \cdots.\]
 That is, let $\overline{\mathbb{T}}=\left(\sum_{i\in \mathbb{Z}} X_{i}\times [0,1]\right)/\sim$
be the quotient of the disjoint union of the spaces $X_i\times[0,1]$
 with $(x,1)\sim (f_{i}(x),0)$ for all $i\in \mathbb{Z}$ and $x\in X_{i}$. Let $\overline{h}: \sum_{i\in \mathbb{Z}} X_{i}\times [0,1] \rightarrow \overline{\mathbb{T}}$ denote the quotient map. We will identify each $X_{i}$ with $\overline{h}(X_{i}\times\{0\})\subseteq \overline{\mathbb{T}}$. As the basepoint for $\overline{\mathbb{T}}$ we take $\ast\in X_0$.

 Let $q:(\overline{\mathbb{T}},\ast) \rightarrow (\mathbb{T},\ast)$ be the unique map making the following diagram commute, where $\left(\sum_{i\in \mathbb{Z}} id \right)|_{X_i\times [0,1]}=id_{\mathbb{E}\times [0,1]}$ for all $i\in\mathbb{Z}$:
\[
\xymatrix@C+2pc{\sum_{i\in \mathbb{Z}} X_{i}\times [0,1]  \ar[r]^{\sum_{i\in \mathbb{Z}} id} \ar[d]^{\overline{h}} &   \mathbb{E}\times[0,1] \ar[d]^{h}\\\overline{\mathbb{T}} \ar[r]^{q} & \mathbb{T}}
\]

Observe that $q:\overline{\mathbb{T}} \rightarrow \mathbb{T}$ is a covering projection. (For every $t\in (0,1)$, the set $U=h(\mathbb{E}\times\left([0,1]\setminus\{t\}\right))$ is open in $\mathbb{T}$, its preimage $q^{-1}(U)$ equals the union of the disjoint open subsets $V_i=\overline{h}((X_{i-1}\times(t,1])\cup (X_i\times[0,t)))$ of $\overline{\mathbb{T}}$ with $i\in \mathbb{Z}$, and $q|_{V_i}:V_i\rightarrow U$ is a homeomorphism for all $i\in \mathbb{Z}$.)

  \begin{figure}[h]
\hspace{-1.2in} \parbox{3in}{\includegraphics[scale=0.5]{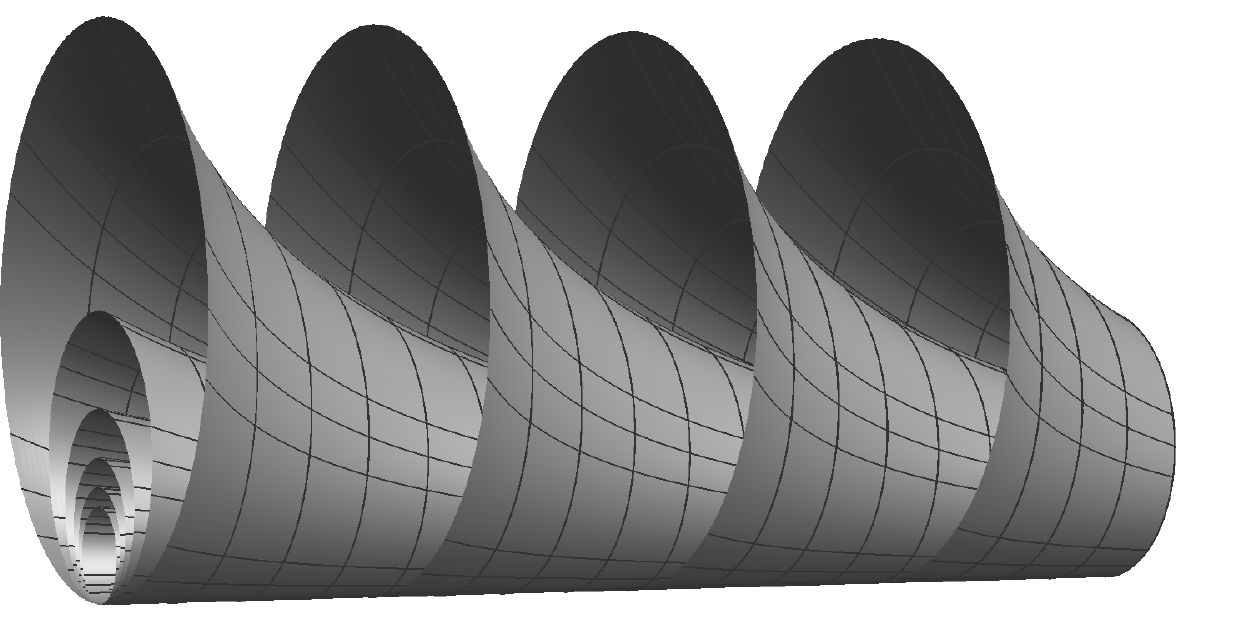}}  \hspace{-.8in} \parbox{1in}{\includegraphics[scale=0.5]{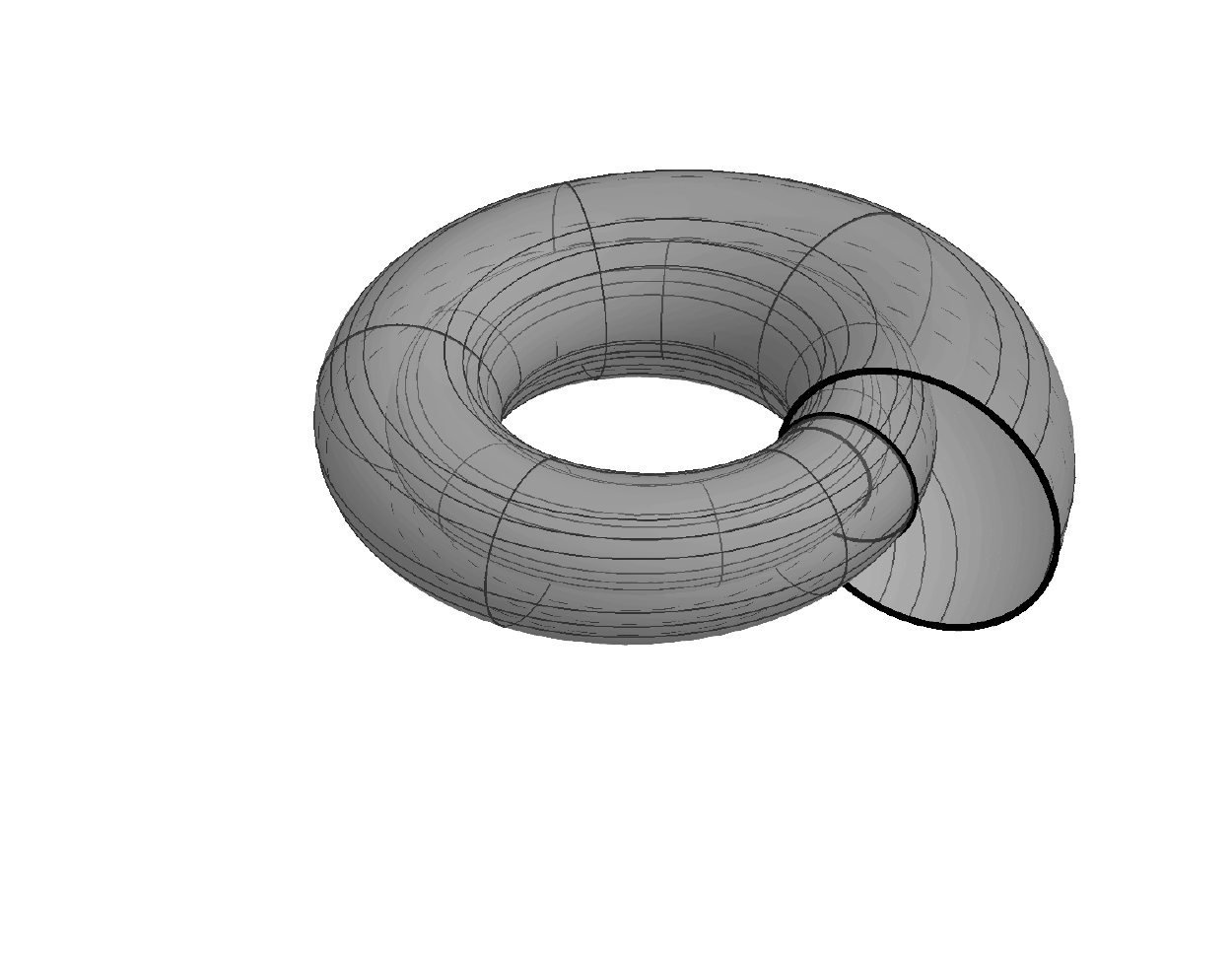}}
\vspace{-.3in}
\caption{A portion of $\overline{\mathbb{T}}$ (left) and  $\mathbb{E}\subseteq \mathbb{T}$ (right).}

 \end{figure}

\pagebreak

\begin{lemma}\label{inj} The homomorphism $\iota_\#:\pi_1(\mathbb{E},\ast)\rightarrow \pi_1(\mathbb{T},\ast)$ is injective.
\end{lemma}

\begin{proof}
Let $[\alpha]\in \pi_1(\mathbb{E},\ast)$ such that $[\alpha]=1\in \pi_1(\mathbb{T},\ast)$. Then $\alpha$ lifts to a loop $\overline{\alpha}:([0,1],\{0,1\})\rightarrow (X_0,\ast)\subseteq (\overline{\mathbb{T}},\ast)$ with $q\circ \overline{\alpha}=\alpha$ and $[\overline{\alpha}]=1\in \pi_1(\overline{\mathbb{T}},\ast)$. Since any homotopy contracting $\overline{\alpha}$ in $\overline{\mathbb{T}}$ has compact image, and since for each $i\in\mathbb{Z}$, $\overline{h}(X_i\times[0,1])\subseteq \overline{\mathbb{T}}$ canonically deformation retracts onto $f_i(X_i)\subseteq X_{i+1}$ along the second coordinate, there is a $k\in \mathbb{N}$ such that $f_k\circ f_{k-1}\circ\cdots\circ f_0\circ\overline{\alpha}$ contracts in $X_{k+1}$. However, $f_\#:\pi_1(\mathbb{E},\ast)\rightarrow \pi_1(\mathbb{E},\ast)$ is injective, because $f:\mathbb{E}\rightarrow \mathbb{E}$ is a section for the quotient map $g:\mathbb{E}\rightarrow \mathbb{E}$ with $g\circ \ell_{m+1}=\ell_m$ for all $m\in \mathbb{N}$ and $g(C_1)=\{\ast\}$; that is, $g\circ f=id_{\mathbb{E}}$.
Hence, $[\overline{\alpha}]=1\in \pi_1(X_0,\ast)$ and $[\alpha]=[q\circ \overline{\alpha}]=1\in \pi_1(\mathbb{E},\ast)$.
\end{proof}

\begin{lemma} \label{SG} We have $\iota_\#\pi_1(\mathbb{E},\ast)\leqslant \pi^s(\mathbb{T},\ast)$. In particular, $\pi^s(\mathbb{T},\ast)$ is uncountable.
\end{lemma}

\begin{proof}
Let $\alpha:([0,1],\{0,1\})\rightarrow (\mathbb{E},\ast)$ be a loop and let $\mathcal U$ be an open cover of $\mathbb{T}$. Choose $U\in {\mathcal U}$ with $\ast\in U$. Choose $k\in\mathbb{N}$ such that $\bigcup_{m=k+1}^\infty C_m\subseteq U$. Then $\delta=f^k\circ\alpha$ lies in $U$ and $[\alpha]=[\lambda]^k[\delta][\lambda]^{-k}\in \pi({\mathcal U},\ast)$. Hence $[\alpha]\in \pi^s(\mathbb{T},\ast)$.

 It is well known that the group $\pi_1(\mathbb{E},\ast)$ is uncountable; see, for example, \cite{CC}. Hence, $\pi^s(\mathbb{T},\ast)$ is uncountable by Lemma~\ref{inj}.
\end{proof}

\begin{lemma}\label{conj}
For every $[\alpha]\in \pi_1(\mathbb{T},\ast)$ there are $i>0$, $j<0$, and a loop $\beta:([0,1],\{0,1\})\rightarrow (\mathbb{E},\ast)\subseteq (\mathbb{T},\ast)$ such that $[\alpha]=[\lambda]^i[\beta][\lambda]^j$.
\end{lemma}

\begin{proof}
Let $\overline{\alpha}:([0,1],0)\rightarrow (\overline{\mathbb{T}},\ast)$ be the lift with $q\circ\overline{\alpha}=\alpha$. Choose $k\in\mathbb{Z}$ such that $\overline{\alpha}(1)\in X_k$. Let $\lambda^k$ be the $|k|$-fold concatenation of $\lambda$ (if $k\geqslant 0$) or of $\lambda^-$ (if $k<0$) and let $\overline{\gamma}:([0,1],0)\rightarrow (\overline{\mathbb{T}},\ast)$ be the lift with $q\circ \overline{\gamma}=\lambda^k$. Then $\overline{\alpha}\cdot (\overline{\gamma})^-$ is a loop. As in the proof of Lemma~\ref{inj}, there is an $i\in \mathbb{Z}$ with $i>\max\{0,k\}$ such that the loop $\overline{\alpha}\cdot (\overline{\gamma})^-$ is freely homotopic within $\overline{\mathbb{T}}$ to a loop $\overline{\beta}$ in $X_i$ with $q\circ \overline{\beta}(0)=q\circ \overline{\beta}(1)=\ast$. Put $\beta=q\circ\overline{\beta}$. Then $[\alpha][\lambda]^{-k}=[\lambda]^i[\beta][\lambda]^{-i}$.
\end{proof}

\begin{proposition} The map $q:\overline{\mathbb{T}} \rightarrow \mathbb{T}$ is a universal covering projection and $q_\#\pi_1(\overline{\mathbb{T}},\ast)= \pi^s(\mathbb{T},\ast)$. In particular, $\pi_1(\overline{\mathbb{T}},\ast)$ is uncountable.
\end{proposition}

\begin{proof} It follows from  the proof of Lemma~\ref{conj} (with $k=0$ and $j=-i$) and Lemma~\ref{SG} that $q_\#\pi_1(\overline{\mathbb{T}},\ast)\subseteq \pi^s(\mathbb{T},\ast)$. In turn, recall from the introduction that for every covering projection $r:(E,e)\rightarrow (\mathbb{T},\ast)$, we have $\pi^s(\mathbb{T},\ast)\subseteq r_\#\pi_1(E,e)$, so that $q_\#\pi_1(\overline{\mathbb{T}},\ast)= \pi^s(\mathbb{T},\ast)\subseteq r_\#\pi_1(E,e)$. Hence, for every covering projection $r:(E,e)\rightarrow (\mathbb{T},\ast)$, there is
a  lift (and covering projection) $s:(\overline{\mathbb{T}},\ast)\rightarrow (E,e)$ with $r\circ s=q$. It follows from Lemma~\ref{SG} that $\pi_1(\overline{\mathbb{T}},\ast)$ is uncountable.
\end{proof}

\section{Covering projections $q_n:\overline{\mathbb{T}}_n\rightarrow \overline{\mathbb{T}}$ ($n\in \mathbb{N}$)}

\noindent
For each $i\in \mathbb{Z}$, there is an embedded infinite ``funnel'' $F_i\subseteq \overline{\mathbb{T}}$,  homeomorphic to $[0,\infty)\times S^1$ and starting at $C_1\subseteq X_i$, with  $\overline{\mathbb{T}}=\bigcup_{i\in \mathbb{Z}} F_i$. Specifically, $F_i$  is the one-sided mapping telescope of the restrictions \[ C_1\stackrel{f_i|_{C_1}}{\longrightarrow} C_2\stackrel{f_{i+1}|_{C_2}}{\longrightarrow} C_3\stackrel{f_{i+2}|_{C_3}}{\rightarrow} \cdots.\]

For each fixed $n\in \mathbb{N}$, we wish to create a covering projection\linebreak $q_n:(\overline{\mathbb{T}}_n,\ast)\rightarrow (\overline{\mathbb{T}},\ast)$ that unwraps the $n$ successive funnels $F_1, F_2, \dots, F_n$ of $\overline{\mathbb{T}}$ into copies of half-planes $[0,\infty)\times \mathbb{R}$. Unwrapping $F_1$ will create infinitely many copies of $\operatorname{cl}\left(\overline{\mathbb{T}}\setminus F_1\right)$, each attached to the resulting half-plane along a horizontal ray. Then there are infinitely many copies of $F_2$ to unwrap, and so on.

\vspace{-.1in}

\begin{figure}[h]
\includegraphics[scale=0.8]{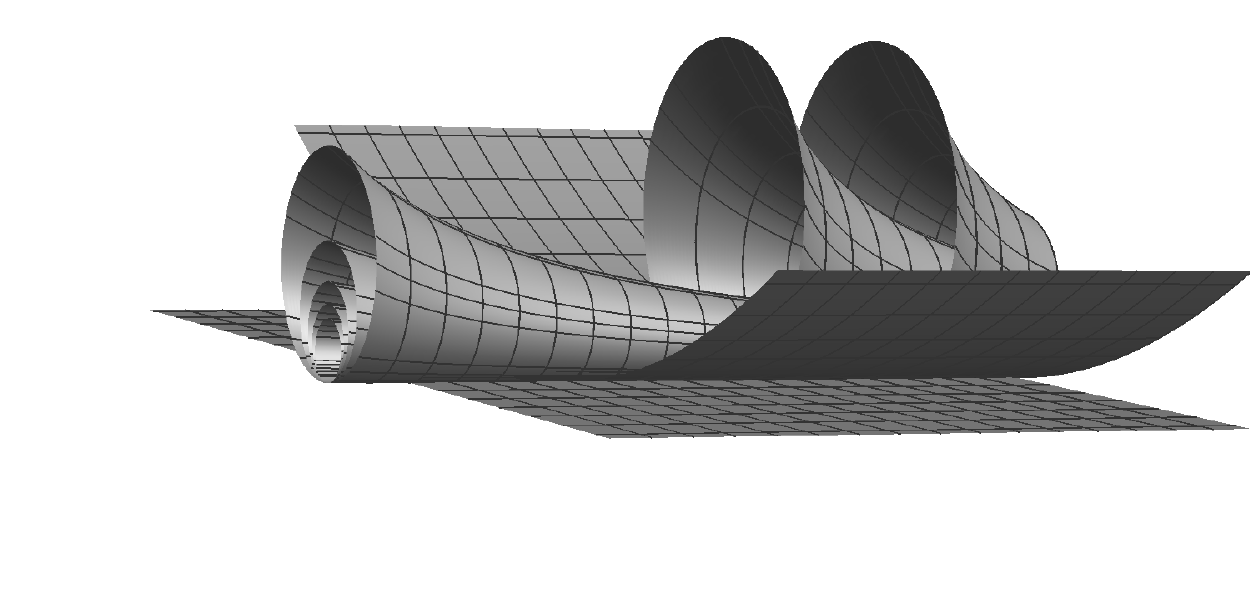}\vspace{-.5in}

\caption{\small Unwrapping the funnels $F_1$ and $F_2$ (local view).}
\end{figure}

 In order to facilitate this, we first prepare the necessary covering projections $p_{n,i}:X_{n,i}\rightarrow X_i$.

 In the three steps below, we will define for each $i\in \mathbb{Z}$, a pointed space $(X_{n,i},\ast)$, a covering projection $p_{n,i}:(X_{n,i},\ast)\rightarrow (X_i,\ast)$, and a map $f_{n,i}:(X_{n,i},\ast)\rightarrow (X_{n,i+1},\ast)$,  such that the following diagram commutes:
\renewcommand{\theequation}{${\mathcal E}_n$}
\begin{equation}\label{commute}
\xymatrix{\cdots \ar[r]^{f_{n,-2}}&X_{n,-1} \ar[d]^{p_{n,-1}} \ar[r]^{f_{n,-1}}  &X_{n,0} \ar[d]^{p_{n,0}} \ar[r]^{f_{n,0}} & X_{n,1} \ar[d]^{p_{n,1}} \ar[r]^{f_{n,1}}  & X_{n,2} \ar[d]^{p_{n,2}} \ar[r]^{f_{n,2}} & \cdots \\
\cdots  \ar[r]^{f_{-2}}&X_{-1} \ar[r]^{f_{-1}} & X_{0}  \ar[r]^{f_{0}} &X_1 \ar[r]^{f_1} & X_{2}  \ar[r]^{f_{2}}  &\cdots}
\end{equation}

First, for $j\geqslant i \geqslant 1$, let $q_{i,j}:(T_{i,j},\ast)\rightarrow (C_i\cup C_{i+1}\cup \cdots \cup C_j,\ast)$ be the universal covering projection and let $\overline{\mathbb{E}}_{i,j}$ denote the tree $T_{i,j}$ with  one copy of $C_1\cup C_2\cup \cdots \cup C_{i-1}\cup C_{j+1}\cup C_{j+2}\cup \cdots$ attached to each vertex along the basepoint. Then $q_{i,j}$ extends canonically to a covering projection $q_{i,j}:(\overline{\mathbb{E}}_{i,j},\ast) \rightarrow (\mathbb{E},\ast)$.
\vspace{5pt}

{\bf Step 1:} $i\geqslant n$. Put
$X_{n,i}=\overline{\mathbb{E}}_{i-n+1,i}$ and  $p_{n,i}=q_{i-n+1,i}:X_{n,i}\rightarrow X_i$. Let $f_{n,i}:T_{i-n+1,i}\rightarrow T_{i-n+2,i+1}$ be the unique homeomorphism such that the following diagram commutes:
\[\xymatrix{(T_{i-n+1,i},\ast) \ar[r]^{f_{n,i}} \ar[d]^{p_{n,i}|_{T_{i-n+1,i}}} & (T_{i-n+2,i+1},\ast) \ar[d]^{p_{n,i+1}|_{T_{i-n+2,i+1} }} \\
(X_i,\ast) \ar[r]^{f_i} & (X_{i+1},\ast)}\]
Let $L_{n,i}=\mathbb{N}\setminus\{i-n+1, i-n+2, \cdots, i\}$. For every $x\in p_{n,i}^{-1}(\ast)$ and for every $m\in \mathbb{N}$, the lift  $\ell_{x,m}:([0,1],0)\rightarrow (X_{n,i},x)$ of the loop $\ell_m:([0,1],\{0,1\})\rightarrow (X_i,\ast)$ with $ p_{n,i}\circ \ell_{x,m}=\ell_m$ is a loop if and only if $m\in L_{n,i}$. Moreover,  for  $m\in \mathbb{N}$, we have $m+1\in L_{n,i+1} \Leftrightarrow m\in L_{n,i}$. Therefore, $f_{n,i}$ extends uniquely to the desired map $f_{n,i}:X_{n,i}\rightarrow X_{n,i+1}$.

\vspace{5pt}

{\bf Step 2:}  $0\leqslant i<n$. We define $X_{n,i}$ by decreasing induction: assuming $X_{n,i+1}$ and $p_{n,i+1}:X_{n,i+1}\rightarrow X_{i+1}$ have been defined, let $\mathcal{X}_{n,i+1}$ denote the set of components of the subset $p_{n,i+1}^{-1}(C_2\cup C_3\cup \cdots)\subseteq X_{n,i+1}$. For each $D\in {\mathcal X}_{n,i+1}$, let $X_{n,i,D}$ be a copy of the tree $T_{1,i}$ with one copy of $C_{i+1}\cup C_{i+2}\cup \cdots$ attached to each vertex along the basepoint. We then define $X_{n,i}=\sum_{D\in {\mathcal X}_{n,i+1}}X_{n,i,D}$ and let $p_{n,i}:X_{n,i}\rightarrow X_i$ be the canonical covering projection. As in Step 1, but with $L_{n,i}=\mathbb{N}\setminus\{1,2, \cdots, i\}$ and  $L_{n,i+1}=\mathbb{N}\setminus\{1,2, \cdots, i+1\}$, we may choose a lift $f_{n,i}$ such that diagram~(\ref{commute}) commutes and  $f(X_{n,i,D})=D$ for all $D\in {\mathcal X}_{n,i+1}$. (Here, $L_{i,0}=\mathbb{N}$.) Put $\ast=f^{-1}_{n,i}(\ast)$.

\vspace{5pt}
{\bf Step 3:}  $i<0$. Note that each component of $X_{n,0}$ is a copy of $\mathbb{E}$. Let $(X_{n,i},\ast)=(X_{n,0},\ast)$, and $p_{n,i}=p_{n,0}$. This time, $L_{n,i}=L_{n,i+1}=\mathbb{N}$. For each component $D$ of $X_{n,i}$ define $f_{n,i}|_D=f_i$.
\vspace{10pt}

\begin{figure}[h]
\includegraphics[scale=0.8]{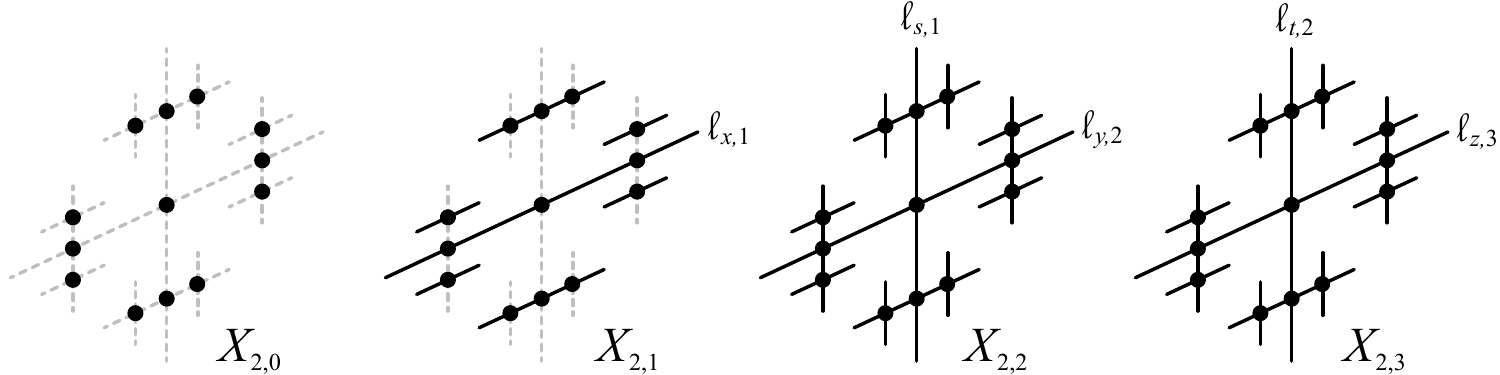}
\caption{\small The non-loops of $X_{n,i}$ for $n=2$ and $i=0,1,2,3$.
 $X_{2,0}$: none;\;\; $X_{2,1}$: $\ell_{x,1}$ with $x\in p_{2,1}^{-1}(\ast)$;\;\; $X_{2,2}$: $\ell_{s,1}, \ell_{y,2}$ with $s\in p_{2,2}^{-1}(\ast)$ and $y=f_{2,1}(x)$;\;\; $X_{2,3}$: $\ell_{t,2}, \ell_{z,3}$ with $t=f_{2,2}(s)$ and $z=f_{2,2}(y)$.}
\end{figure}

Finally,  let $\overline{\mathbb{T}}_n$ be the mapping telescope of the bi-infinite sequence
\[\cdots \stackrel{f_{n,-2}}{\longrightarrow} X_{n,-1} \stackrel{f_{n,-1}}{\longrightarrow} X_{n,0} \stackrel{f_{n,0}}{\longrightarrow} X_{n,1} \stackrel{f_{n,1}}{\longrightarrow}  X_{n,2} \stackrel{f_{n,2}}{\longrightarrow} \cdots\]
 with basepoint $\ast\in X_{n,0}$. That is, $\overline{\mathbb{T}}_n=\left(\sum_{i\in \mathbb{Z}} X_{n,i}\times [0,1]\right)/\sim$ with $(x,1)\sim (f_{n,i}(x),0)$ for all $i\in \mathbb{Z}$ and $x\in X_{n,i}$. Note that $\overline{\mathbb{T}}_n$ is connected.
  Let $\overline{h}_n: \sum_{i\in \mathbb{Z}} X_{n,i}\times [0,1] \rightarrow \overline{\mathbb{T}}_n$ denote the quotient map. We will identify each $X_{n,i}$ with $\overline{h}_n(X_{n,i}\times\{0\})\subseteq \overline{\mathbb{T}}_n$.

  Let $q_n:(\overline{\mathbb{T}}_n,\ast) \rightarrow (\mathbb{T},\ast)$ be the unique map making the following diagram commute:
\[
\xymatrix@C+3pc{\sum_{i\in \mathbb{Z}} X_{n,i}\times [0,1]  \ar[r]^{\sum_{i\in \mathbb{Z}} p_{n,i}\times id} \ar[d]^{\overline{h}_n} &  \sum_{i\in \mathbb{Z}} X_{i}\times [0,1] \ar[d]^{\overline{h}}\\\overline{\mathbb{T}}_n \ar[r]^{q_n} & \overline{\mathbb{T}}}
\]
Observe that $q_n:\overline{\mathbb{T}}_n \rightarrow \overline{\mathbb{T}}$ is a covering projection.

\section{Main theorem}

The key property of $\mathbb{T}$ is that it contains loops requiring a composition of at least two covering projections to unwind them:

\begin{theorem}\label{main}$ $

{\em (a)} For every loop $\alpha:([0,1],\{0,1\})\rightarrow (\mathbb{E},\ast)\subseteq (\mathbb{T},\ast)$ and every covering projection $p:(E,e)\rightarrow (\mathbb{T},\ast)$, the  lift $\alpha':([0,1],0)\rightarrow (E,e)$ with $p\circ\alpha'=\alpha$ is a loop, i.e. $\alpha'(1)= e$.

{\em (b)} For every essential loop $\alpha:([0,1],\{0,1\})\rightarrow (\mathbb{T},\ast)$ there are two covering projections
$p_1:(E_1,e_1)\rightarrow (\mathbb{T},\ast)$ and $p_2:(E_2,e_2)\rightarrow (E_1,e_1)$ such that the  lift $\alpha'':([0,1],0)\rightarrow (E_2,e_2)$ with  $p_1\circ p_2\circ\alpha''=\alpha$ is not a loop, i.e. $\alpha''(1)\not=e_2$.
\end{theorem}

\begin{proof}
(a) This follows from Lemma~\ref{SG} and the characterization of $\pi^s(\mathbb{T},\ast)$ from the introduction.

(b) By Lemma~\ref{conj}, $[\alpha]=[\lambda]^i[\beta][\lambda]^j$ for some $i>0$, some $j<0$, and some loop $\beta:([0,1],\{0,1\})\rightarrow (\mathbb{E},\ast)\subseteq (\mathbb{T},\ast)$. If $\beta$ is null-homotopic in $\mathbb{E}$ or if $i+j\not=0$, then the lift $\overline{\alpha}:([0,1],0)\rightarrow (\overline{\mathbb{T}},\ast)$ with $q\circ \overline{\alpha}=\alpha$ is not a loop. We may therefore assume that $\beta$ is essential in $\mathbb{E}$ and that $j=-i$.

For each $n\in \mathbb{N}$, let $r_n:\mathbb{E}\rightarrow \bigcup_{m=1}^n C_m$ denote the retraction with $r_n(\bigcup_{m=n+1}^\infty C_m)=\{\ast\}$. Then there is an $n\in \mathbb{N}$ such that $r_n\circ \beta$ is essential in $\bigcup_{m=1}^n C_m$ \cite{CC}. We may assume that $n>i$.
Since $\pi_1(\mathbb{E},\ast)$ is the free product of $\pi_1(\bigcup_{m=1}^n C_m,\ast)$ and $\pi_1(\bigcup_{m=n+1}^\infty C_m,\ast)$, there are loops $\gamma_1, \gamma_2, \dots, \gamma_s$ in $\bigcup_{m=1}^n C_m$ based at $\ast$ and loops $\delta_1, \delta_2, \dots, \delta_s$ in $\bigcup_{m=n+1}^\infty C_m$ based at $\ast$, such that $[\beta]=[\gamma_1\cdot \delta_1\cdot \gamma_2\cdot \delta_2\cdots \gamma_s\cdot \delta_s]\in \pi_1(\mathbb{E},\ast)$
and  $[r_n\circ \beta]=[\gamma_1\cdot\gamma_2\cdots \gamma_s]\not=1\in \pi_1(\bigcup_{m=1}^n C_m,\ast)$. Let $p_1:(E_1,e_1)\rightarrow (\mathbb{T},\ast)$ be the covering projection  $q:(\overline{\mathbb{T}},\ast)\rightarrow (\mathbb{T},\ast)$. Let $\widehat{\lambda}:([0,1],0)\rightarrow (\overline{\mathbb{T}}_n,\ast)$ be the lift of $\lambda^{n-i}$ with $q \circ q_n\circ \widehat{\lambda}=\lambda^{n-i}$.  Put $E_2=\overline{\mathbb{T}}_n$, $e_2=\widehat{\lambda}(1)$, and $p_2=q_n:(E_2,e_2)\rightarrow (E_1,e_1)$.

We now consider the lift of the loop $\lambda^i\cdot \gamma_1\cdot \delta_1\cdot \gamma_2\cdot \delta_2\cdots \gamma_s\cdot \delta_s\cdot \lambda^{-i}$ to $E_2=\overline{\mathbb{T}}_n$ through both covering projections $p_1:(E_1,e_1)\rightarrow (\mathbb{T},\ast)$ and $p_2:(E_2,e_2)\rightarrow (E_1,e_1)$, starting at $e_2$. The lift of $\lambda^i$ starting at $e_2$ ends at the basepoint of the tree $T_{1,n}\subseteq X_{n,n}$. By construction of $X_{n,n}$ and the choice of our $\gamma's$ and $\delta's$, the subsequent lift of $\gamma_1\cdot \delta_1\cdot \gamma_2\cdot \delta_2\cdots \gamma_s\cdot \delta_s$ starting at the basepoint of $T_{1,n}$ is not a loop, so that the entire lift is not a loop. We conclude that the lift $\alpha'':([0,1],0)\rightarrow (E_2,e_2)$ with  $p_1\circ p_2\circ\alpha''=\alpha$ is not a loop.
\end{proof}

\begin{corollary}\label{fromproof}
For every $1\not=[\overline{\alpha}]\in \pi_1(\overline{\mathbb{T}},\ast)$,
there is a covering projection $r:(E,e)\rightarrow (\overline{\mathbb{T}},\ast)$ such that $[\overline{\alpha}]\not\in r_\#\pi_1(E,e)$.
\end{corollary}

\begin{proof} See proof of Theorem~\ref{main}(b).
\end{proof}

\begin{corollary}
We have $\pi^s(\overline{\mathbb{T}},\ast)=1$, while no universal covering projection over
$\overline{\mathbb{T}}$ exists.
\end{corollary}

\begin{proof}
It follows from Corollary~\ref{fromproof} and the characterization of $\pi^s$ from the introduction that $\pi^s(\overline{\mathbb{T}},\ast)=1$. If there were a universal covering projection $u:(E,e)\rightarrow (\overline{\mathbb{\mathbb{T}}},\ast)$, then $u_\#\pi_1(E,e)=\pi^s(\overline{\mathbb{T}},\ast)=1$ and $\overline{\mathbb{T}}$ would be semilocally simply connected. However, it follows from the proof of Lemma~\ref{inj} that $\overline{\mathbb{T}}$ is not semilocally simply connected.
\end{proof}

\begin{corollary}\label{simply} The universal object in the category of fibrations with unique path lifting and path-connected total space over~$\mathbb{T}$ has simply connected domain, although $\mathbb{T}$ has uncountable Spanier group.
\end{corollary}

\begin{proof}
As discussed in \cite[end of \S2.5]{S}, we can form an indexed  set $\{p_i:(Y_i,y_i)\rightarrow (\mathbb{T},\ast)\mid i\in I\}$
containing one representative for each pointed fibration with unique path lifting and path-connected total space over $(\mathbb{T},\ast)$. Let $\widetilde{\mathbb{T}}$ be the path component of the fibered product $\{(x_i)_{i \in I} \in \prod_{i\in I} Y_i\mid p_i(x_i)=p_j(x_j)$ for all $i, j\in I\}$  with basepoint $(y_i)_{i\in I}\in \widetilde{\mathbb{T}}$ and let $p:\widetilde{\mathbb{T}}\rightarrow \mathbb{T}$ be defined by $p((x_i)_{i \in I})=p_i(x_i)$ for any $i\in I$. Then $p:(\widetilde{\mathbb{T}},\ast)\rightarrow (\mathbb{T},\ast)$ is the universal object in the category of (pointed) fibrations with unique path lifting and path-connected total space over $(\mathbb{T},\ast)$.

We show that $p_\#\pi_1(\widetilde{\mathbb{T}},\ast)=1$.  Suppose, to the contrary, that there exists $1\not=[\alpha]\in p_\#\pi_1(\widetilde{\mathbb{T}},\ast)\leqslant \pi_1(\mathbb{T},\ast)$. By Theorem~\ref{main}, there are covering projections
$p_1:(E_1,e_1)\rightarrow (\mathbb{T},\ast)$ and $p_2:(E_2,e_2)\rightarrow (E_1,e_1)$ such that the lift $\alpha'':([0,1],0)\rightarrow (E_2,e_2)$ with  $p_1\circ p_2\circ\alpha''=\alpha$ is not a loop. Although  $p_1\circ p_2:(E_2,e_2)\rightarrow (\mathbb{T},\ast)$ may not be a covering projection, it is a fibration with unique path lifting. Therefore, there is a map $p':(\widetilde{\mathbb{T}},\ast)\rightarrow (E_2,e_2)$ such that $p_1 \circ p_2 \circ p'=p$. Since $\alpha$ does not lift to a loop of $E_2$ at $e_2$, it does not lift to a loop of $\widetilde{\mathbb{T}}$ at $\ast$, contradicting the fact that $[\alpha]\in p_\#\pi_1(\widetilde{\mathbb{T}},\ast)$.
\end{proof}

\begin{corollary}\label{notInvLim}
The universal object in the category of fibrations with unique path lifting and path-connected total space over~$\mathbb{T}$ is not a path component of any inverse limit of covering projections over $\mathbb{T}$.
\end{corollary}

\begin{proof}
Let $(\Lambda,\leqslant)$ be a directed set. For $i\in \Lambda$, let $p_i:(Y_i,y_i)\rightarrow (\mathbb{T},\ast)$ be a covering projection with connected total space $Y_i$. For each pair $i,j\in \Lambda$, with $i\leqslant j$ let $f_{i,j}:(Y_j,y_j)\rightarrow (Y_i,y_i)$ be a (covering projection) map such that $p_i\circ f_{i,j}=p_j$ and such that $f_{i,j}\circ f_{j,k}=f_{i,k}$ whenever $i\leqslant j\leqslant k$.
Let $u:(Y,y)\rightarrow (\mathbb{T},\ast)$ be the inverse limit of the inverse system $\left(p_i, f_{i,j},\Lambda\right)$ in the category of pointed fibrations with unique path lifting over $(\mathbb{T},\ast)$, and note that $ u_\#\pi_1(Y,y)=\bigcap_{i\in \Lambda} (p_i)_\#\pi_1(Y_i,y_i)$.
Since $1\not=\pi^s(\mathbb{T},\ast)\leqslant \bigcap_{i\in \Lambda} (p_i)_\#\pi_1(Y_i,y_i)$, the path component of $y$ in $Y$ is not simply connected. The claim now follows from Corollary~\ref{simply}.
\end{proof}

\section{Some final remarks}
\renewcommand{\theequation}{$\mathcal L$}

\noindent We recall from \cite[Theorem 2.4.5]{S} that every fibration $u:Y\rightarrow X$ with unique path lifting satisfies the following lifting property:
\begin{equation} \label{LiftingProperty}
\parbox{4.2in}{\em
For every connected and locally path-connected space $Z$, for every $y\in Y$, and for every continuous function \linebreak $f:(Z,z)\rightarrow (X,u(y))$ with $f_{\#}\pi_1(Z,z)\leqslant u_{\#}\pi_1(Y,y)$, there is a unique continuous lift $f':(Z,z)\rightarrow (Y,y)$ with $u\circ f'=f$.
}
\end{equation}
Note that a continuous function $u:(Y,\ast)\rightarrow (X,\ast)$ with connected and locally path-connected $Y$ which satisfies lifting property~(\ref{LiftingProperty}) is uniquely determined by the subgroup $u_{\#}\pi_1(Y,\ast)$ of $\pi_1(X,\ast)$.

\begin{remark} \label{nonunique} For every connected and locally path-connected space $X$, and $x\in X$, there is a fibration $u:(Y,y)\rightarrow (X,x)$ with unique path lifting and path-connected $Y$ such that $u_\#\pi_1(Y,y)=\pi^s(X,x)$. However, such a fibration need not be unique, since $Y$ need not be locally path connected. (For example, for every $\mathcal U\in Cov(X)$, we may choose a regular covering projection $q:(E,e)\rightarrow (X,x)$ with connected total space $E$ such that $q_\#\pi_1(E,e)=\pi({\mathcal U},x)$ \cite[2.5.13]{S}. Lifting between these covering projections whenever possible, creates an inverse system of pointed regular covering projections over $(X,x)$ which is directed inversely to inclusion of the corresponding subgroups $\pi({\mathcal U},x)$.  Let $u:(Y',y)\rightarrow (X,x)$ be the inverse limit as in Corollary~\ref{notInvLim} and let $Y$ be the path component of $Y'$ with $y\in Y$.)
\end{remark}

\begin{remark}\label{GenUnivCov} In cases where a simply connected universal fibration with unique path lifting does not exist, there might still be other universal lifting objects that are unique and simply connected.

Indeed, as discussed in \cite[Example 6.10]{BFi}, there is a Peano continuum $B$ with the following properties: On the one hand, there is a (unique)
map  $g:(\widehat{B},\widehat{b})\rightarrow (B,b)$ with
connected and locally path-connected $\widehat{B}$ that satisfies lifting property (\ref{LiftingProperty}) and
 has  $g_\#\pi_1(\widehat{B},\widehat{b})=1$. On the other hand, the singleton $\{1\}$ is not closed in $\pi_1(B,b)$, when $\pi_1(B,b)$ is given the quotient topology of the loop space $\Omega(B,b)$ with the compact-open topology under the canonical map $[\,\cdot\, ]:\Omega(B,b)\rightarrow \pi_1(B,b)$.

If $v:(Z,z)\rightarrow (B,b)$ is any fibration with unique path lifting, then
the fiber $v^{-1}(b)$ does not contain nonconstant paths and is therefore $T_1$.
Consequently, since fibrations have continuous lifting of paths \cite[2.7.8]{S}, $v_\#\pi_1(Z,z)$ is closed in $\pi_1(B,b)$. So, if $p:(\widetilde{B},\widetilde{b})\rightarrow (B,b)$ denotes the universal object in the category of fibrations with unique path lifting and path-connected total space over $B$, then $1=g_\#\pi_1(\widehat{B},\widehat{b}) \lvertneqq p_\#\pi_1(\widetilde{B},\widetilde{b})\leqslant \pi^s(B,b)$.
\end{remark}

\begin{remark}
In \cite{FZ2007} it is shown that for every path-connected space $X$, there is a (unique) map $g:(\widehat{X},\widehat{x})\rightarrow (X,x)$ with connected and locally path-connected  $\widehat{X}$ that
satisfies lifting property (\ref{LiftingProperty}) and has $g_\#\pi_1(\widehat{X},\widehat{x})=\pi^s(X,x)$. However, in general, $g:(\widehat{X},\widehat{x})\rightarrow (X,x)$ is neither a covering projection nor a fibration.
\end{remark}

\end{document}